\documentclass[8 pt]{amsproc}
\usepackage{amssymb}
\usepackage{blindtext, rotating}
\usepackage[english]{babel}
 
\oddsidemargin=-1cm\evensidemargin=-1cm
\begin{document}
\numberwithin{equation}{section}

\def\1#1{\overline{#1}}
\def\2#1{\widetilde{#1}}
\def\3#1{\widehat{#1}}
\def\4#1{\mathbb{#1}}
\def\5#1{\frak{#1}}
\def\6#1{{\mathcal{#1}}}

\def\C{{\4C}}
\def\R{{\4R}}
\def\N{{\4N}}
\def\Z{{\4Z}}

\author{Valentin Burcea }
\title{Real Submanifolds in Complex Spaces:  Upgrades}
 \begin{abstract}  Provided the coordinates $\left(z_{11},z_{12},\dots,z_{1N},\dots,z_{m1},z_{m2},\dots,z_{mN},w_{11},w_{12},\dots,w_{1m},\dots,w_{m1},w_{m2},\dots,w_{mm}\right)$  in $\mathbb{C}^{mN+m^{2}}$  defining   the matrices $W=\left\{w_{ij}\right\}_{1\leq i,j\leq m}$ and $Z=\left\{z_{ij}\right\}_{1\leq i\leq m \atop{1\leq j \leq N}}$, we consider   the Special Class of the Real-Analytic Submanifolds $M$ defined near $p=0$ by  $W=Z\overline{Z}^{t}+\rm{O}\left(3\right)$. We prove that $M$ is biholomorphically equivalent to the  Model
   $W=Z\overline{Z}^{t}$    if and only  if $M$ is formally equivalent to the  Model
   $W=Z\overline{Z}^{t}$. It follows from a non-equidimensional version of this result.
   
  \end{abstract}

\address{V. Burcea: INDEPENDENT}
\email{vdburcea@gmail.com}
\thanks{\emph{Keywords:} Real Submanifold, Equivalence Problem, Fischer Decomposition}
\thanks{Special Thanks  to CAPES at Federal University of Santa Catarina, Brazil;}
\thanks{Emphasizing that the reference \cite{V1} was fully supported by  Science Foundation Ireland, Grant 06/RFP/MAT 018.}
 \maketitle
\def\Label#1{\label{#1}{\bf (#1)}~}

\def\cn{{\C^n}}
\def\cnn{{\C^{n'}}}
\def\ocn{\2{\C^n}}
\def\ocnn{\2{\C^{n'}}}


\def\dist{{\rm dist}}
\def\const{{\rm const}}
\def\rk{{\rm rank\,}}
\def\id{{\sf id}}
\def\tr{{\bf tr\,}}
\def\aut{{\sf aut}}
\def\Aut{{\sf Aut}}
\def\CR{{\rm CR}}
\def\GL{{\sf GL}}
\def\Re{{\sf Re}\,}
\def\Im{{\sf Im}\,}
\def\span{\text{\rm span}}
\def\Diff{{\sf Diff}}

\def\codim{{\rm codim}}
\def\crd{\dim_{{\rm CR}}}
\def\crc{{\rm codim_{CR}}}

\def\phi{\varphi}
\def\eps{\varepsilon}
\def\d{\partial}
\def\a{\alpha}
\def\b{\beta}
\def\g{\gamma}
\def\G{\Gamma}
\def\D{\Delta}
\def\Om{\Omega}
\def\k{\kappa}
\def\l{\lambda}
\def\L{\Lambda}
\def\z{{\bar z}}
\def\w{{\bar w}}
\def\Z{{\1Z}}
\def\t{\tau}
\def\th{\theta}

\emergencystretch15pt \frenchspacing

\newtheorem{Thm}{Theorem}[section]
\newtheorem{Cor}[Thm]{Corollary}
\newtheorem{Pro}[Thm]{Proposition}
\newtheorem{Lem}[Thm]{Lemma}

\theoremstyle{definition}\newtheorem{Def}[Thm]{Definition}

\theoremstyle{remark}
\newtheorem{Rem}[Thm]{Remark}
\newtheorem{Exa}[Thm]{Example}
\newtheorem{Exs}[Thm]{Examples}

\def\bl{\begin{Lem}}
\def\el{\end{Lem}}
\def\bp{\begin{Pro}}
\def\ep{\end{Pro}}
\def\bt{\begin{Thm}}
\def\et{\end{Thm}}
\def\bc{\begin{Cor}}
\def\ec{\end{Cor}}
\def\bd{\begin{Def}}
\def\ed{\end{Def}}
\def\br{\begin{Rem}}
\def\er{\end{Rem}}
\def\be{\begin{Exa}}
\def\ee{\end{Exa}}
\def\bpf{\begin{proof}}
\def\epf{\end{proof}}
\def\ben{\begin{enumerate}}
\def\een{\end{enumerate}}
\def\beq{\begin{equation}}
\def\eeq{\end{equation}}

\section{Introduction and Main Result}

Through this paper, we study Classes of Real-Formal Submanifolds   derived from Shilov Boundaries of Bounded and Symmetric Domains of first kind (see \cite{KiZa1},\cite{KiZa2})  in the light of the   Theorem of Moser\cite{Mo} and its generalizations from   Gong\cite{Go1} and Huang-Yin\cite{HuYi}. These Real Submanifolds are ,,modelled'' by Shilov Boundaries of Bounded and Symmetric Domains of first kind\cite{KaZa1},\cite{KaZa2}. Such Real Submanifolds are named $\mathcal{BSD}$-Manifolds through this paper. They are denoted by 
\begin{equation}  \mathcal{M}_{m,N}:\hspace{0.1 cm}  W=Z\overline{Z}^{t}+\rm{O}\left(3\right)\subset \mathbb{C}^{mN+m^{2}}.  \label{MAN}\end{equation}

In the light of   standard linear embeddings, we consider  coordinates defined by the standard embeddings
\begin{equation} \mathbb{C}^{mN+m^{2}}  \rightarrow  \mathbb{C}^{{m}\left(N+1\right)+{m}^{2}}.  \label{incl}\end{equation}

In the light of (\ref{incl}), we work with the following matrices
\begin{equation}\quad\quad\quad\quad\left(W, Z\right)=\left(\left\{w_{ij}\right\}_{1\leq i,j\leq m},  \left\{z_{ij}\right\}_{1\leq i\leq m\atop{1\leq j \leq N}}\right),\hspace{0.1 cm}\mbox{defining the  $\mathcal{BSD}$-Model} \hspace{0.1 cm}\mathcal{BSD}_{m,N}:\hspace{0.1 cm}  W=Z\overline{Z}^{t}\subset \mathbb{C}^{Nm+m^{2}},\label{coordA}
\end{equation} 

\begin{equation}  \left(W', Z'\right)=\left(\left\{w'_{i'j'}\right\}_{1\leq i',j'\leq m},  \left\{z'_{i'j'}\right\}_{1\leq i'\leq m\atop{1\leq j' \leq N+1}}\right),\hspace{0.1 cm}\mbox{defining the  $\mathcal{BSD}$-Model}\hspace{0.1 cm}\mathcal{BSD}_{m,N+1}:\hspace{0.1 cm}  W'=Z'\overline{{Z'}}^{t}\subset \mathbb{C}^{{m}\left(N+1\right)+{m}^{2}}.  \label{coordC}
\end{equation} 

These settings are used in order to study the   non-equidimensional  version of   the   Equivalence Problem between two 
Real-Analytic Submanifolds in Complex Spaces.  Going back to Poincare\cite{Poincare},  it   asks if two Real-Analytic formally equivalent Submanifolds in Complex Spaces, are  actually biholomorphically equivalent. It is one  of the most beautiful problems in Complex Analysis.  Chern-Moser\cite{CM}  proved the convergence of any Formal   Equivalence   between two Levi-NonDegenerate Real-Analytic Hypersurfaces. In contrast to the results of  Mir\cite{Mir1},\cite{Mir2}   in the CR finite type Setting,   Kossovskiy-Shafikov\cite{RI1},\cite{RI2} proved the existence of Real-Analytic hypersurfaces  formally, but not holomorphically equivalent in the infinite type Setting (see also \cite{KL}).  In the CR Singular Setting (see \cite{bu}),   Moser-Webster\cite{MoWe} and Gong\cite{Go2}  constructed  Real-Analytic Submanifolds   formally equivalent, but not  holomorphically equivalent.

 Formal expansions of Formal (Holomorphic) Mappings are considered, on  entries in the defining equations of matrix type, in order to develop formal constructions of normal form type. It follows that

\bt\label{teo1} Let $M\subset\mathbb{C}^{mN+m^{2}}$ be the Real-Analytic Submanifold defined near $p=0$ by (\ref{MAN}), and the Model (\ref{coordC}). Then $M$ is holomorphically embeddable into the  Model (\ref{coordC}) if and only  if $M$ is formally embeddable into the  Model (\ref{coordC}).
\et

 We  recognize Generalized Fischer Decompositions\cite{Sh}  on homogeneous terms     in order to extend and  reconstruct the roots and the ruins from \cite{V3}. We define iteratively Spaces of Fischer-Normalizations as in \cite{V2}   in order to compute the Formal (Holomorphic) Embedding, excepting  an infinite number of components.  In order to obtain its Convergence, we use  formal automorphisms group of the target $\mathcal{BSD}$-Model in order to cancel by composition its undetermined components.   We reformulate the arguments of   rapid convergence of Moser\cite{Mo} and Huang-Yin\cite{HuYi} in order to solve the non-equidimensional convergence problem.  In particular, we obtain:

\bc\label{c1} Let $M\subset\mathbb{C}^{mN+m^{2}}$ be the Real-Analytic Submanifold defined near $p=0$ by (\ref{MAN}), and the Model
 \begin{equation} W=Z\overline{Z}^{t}.\label{15}\end{equation}
 
Then $M$ is holomorphically equivalent to the  Model (\ref{15}) if and only  if $M$ is formally equivalent to the  Model (\ref{15}).
\ec

This should be the correct statement of Theorem 1.1 from \cite{V3}. It was not entirely proved in \cite{V3}, because the proof of Lemma $2.4$ was just wrong  in \cite{V3}.  Its  proof requires the  implementation of the  strategy from \cite{V4} in order to obtain the desired normalizing automorphism. It must  uniquely determined contrary to the wrong claims from \cite{V3} caused by severe unbalances.

\section{Settings} We work in  the coordinates  
 \begin{equation}\left(w_{11},w_{12},\dots,w_{1m},\dots,w_{m1},w_{m2},\dots,w_{mm};z_{11},z_{12},\dots,z_{1N},\dots,z_{m1},z_{m2},\dots,z_{mN}\right)\in \mathbb{C}^{mN+m^{2}}.\label{coord}\end{equation}

 In particular, we  consider matrices according to  the identifications
\begin{equation}\begin{split}& W=\begin{pmatrix}  w_{11} & w_{12} &\dots &w_{1m} \\ w_{21} & w_{22} &\dots & w_{2m} \\ \vdots &\vdots &\ddots &\vdots  \\ w_{m1} & w_{m2} &\dots & w_{mm}\end{pmatrix}\equiv\left(w_{11},w_{12},  \dots, w_{1m},w_{21},w_{22},  \dots, w_{2m},\dots\dots,w_{m1},w_{m2},  \dots, w_{mm}\right),\\&    Z=\begin{pmatrix} z_{11} &z_{12}& \dots &z_{1N} \\ z_{21} &z_{22}& \dots &z_{2N} \\  \vdots &\vdots &\ddots &\vdots \\ z_{m1} &z_{m2}&\dots &z_{mN}\end{pmatrix} \equiv\left(z_{11}, z_{12}, \dots,z_{1N},z_{21}, z_{22}, \dots,z_{2N},\dots\dots, z_{m1}, z_{m2}, \dots,z_{mN}\right).\end{split}\label{matrix}\end{equation}

In the  regards of (\ref{coord}) and (\ref{matrix}), we use  the identifications
\begin{equation} \begin{split}\left\{1,2,3,\dots,m^{2} \right\}\equiv & \left\{ (1,1),(1,2),\dots,\left(1,m\right)\right.,\\&  \hspace{0.2 cm} (2,1),(2,2),\dots,\left(2,m\right),
\\&\quad\hspace{0.2 cm}\vdots\hspace{0.1 cm}\quad\quad\vdots\quad\hspace{0.2 cm}\ddots\quad\hspace{0.1 cm}\vdots\\& \left.\hspace{0.2 cm}(m,1),(m,2),\dots,\left(m,m\right)\right\},\end{split}\quad\quad\quad \begin{split}\left\{1,2,3,\dots,mN \right\}\equiv & \left\{ (1,1),(1,2),\dots,\left(1,N\right)\right.,\\&  \hspace{0.2 cm} (2,1),(2,2),\dots,\left(2,N\right),
\\&\quad\hspace{0.2 cm}\vdots\hspace{0.1 cm}\quad\quad\vdots\quad\hspace{0.2 cm}\ddots\quad\hspace{0.1 cm}\vdots\\& \left.\hspace{0.2 cm}(m,1),(m,2),\dots,\left(m,N\right)\right\}.\end{split}\label{Ident}
\end{equation}

In particular, we  use the     identifications
\begin{equation}\begin{split}& J:=\begin{pmatrix} j_{11} & j_{12} &\dots &j_{1m} \\ j_{21} & j_{22} &\dots &j_{2m} \\ \vdots & \vdots &\ddots &\vdots  \\ j_{m1} &j_{m2} &\dots &j_{mm}\end{pmatrix}\equiv\left(j_{11},j_{12},  \dots, j_{1m},j_{21},j_{22},  \dots, j_{2m},\dots\dots,j_{m1},j_{m2},  \dots, j_{mm}\right)\in\mathbb{N}^{m^{2}},\\&   I:=\begin{pmatrix} i_{11} &i_{12}& \dots&i_{1N} \\ i_{21} &i_{22}& \dots&i_{2N} \\\vdots &\vdots &\ddots &\vdots \\ i_{m1} &i_{m2}&\dots&i_{mN}\end{pmatrix} \equiv\left(i_{11},i_{12},  \dots,i_{1N},i_{21},i_{22},  \dots,i_{2N},
\dots\dots,i_{m1},i_{m2},  \dots,i_{mN}\right)\in\mathbb{N}^{mN}.\end{split}\label{Ident1}\end{equation}

We define lengths of multi-indexes:
\begin{equation}\begin{split}&\left|J\right|=j_{11}+j_{12}+\dots+j_{1m}+j_{21}+j_{22}+\dots+ j_{2m}+\dots+j_{m1}+j_{m2}+\dots+j_{mm}, \quad\quad\mbox{for all $J\in\mathbb{N}^{m^{2}}$,}\\&\hspace{0.03 cm}
\left|I\right|=i_{11}+i_{12}+  \dots+i_{1N}+i_{21}+i_{22}+\dots+i_{2N}+
\dots\dots+i_{m1}+i_{m2}+\dots+i_{mN}, \quad\hspace{0.1 cm}\mbox{for all $I\in\mathbb{N}^{mN}$.} \end{split}\label{IJ}
\end{equation}

In the light of (\ref{matrix}) and (\ref{Ident1}),  we write  
\begin{equation}\begin{split}&
W^{J}=w_{11}^{j_{11}}w_{12}^{j_{12}}\dots w_{1m}^{j_{1m}} w_{21}^{j_{21}}w_{22}^{j_{22}}\dots w_{2m}^{j_{2m}}\dots\dots w_{m1}^{j_{m1}}w_{m2}^{j_{m2}}\dots w_{mm}^{j_{mm}},
\\&  \hspace{0.14 cm} Z^{I}=z_{11}^{i_{11}}z_{12}^{i_{12}}\dots z_{1N}^{i_{1N}} z_{21}^{i_{21}}z_{22}^{i_{22}}\dots z_{2N}^{i_{2N}}\dots\dots z_{m1}^{i_{m1}}z_{q2}^{i_{m2}}\dots z_{mN}^{i_{mN}}.\end{split}
\label{IJ1}\end{equation}

Regardless of the considered natural numbers, we generalize  the standard hermitian inner-product:
\begin{equation}\left<L, V\right>=L\overline{V}^{t},\quad\mbox{for $L \in \mathcal{M}_{m,n}\left(\mathbb{C}\right)$ and $V \in \mathcal{M}_{n,p}\left(\mathbb{C}\right)$, for all $m,n,p\in\mathbb{N}^{\star}$.}\label{vb} 
\end{equation}

Analogous notations and identifications to (\ref{matrix}),(\ref{Ident}),(\ref{Ident1}),(\ref{IJ}) and (\ref{IJ1})   are considered   in order to work in  the folllowing coordinates  
 \begin{equation}\left({w'}_{11}, \dots,{w'}_{1\hspace{0.05 cm}m},\dots,{w'}_{1\hspace{0.05 cm}m}, \dots,{w'}_{m\hspace{0.05 cm}m};{z'}_{11}, \dots,{z'}_{1\hspace{0.05 cm}N+1},\dots,{z'}_{m\hspace{0.05 cm}1}, \dots,{z'}_{m\hspace{0.05 cm}N+1}\right)\in \mathbb{C}^{m\left(N+1\right)+m^{2}}.\label{coord2}\end{equation}

In  the light of $(2.4)$ from the root\cite{V3}, we use   the   row-vectors

 \begin{equation}  \begin{pmatrix} \mathcal{L}_{1} \\ \mathcal{L}_{2}\\ \vdots \\ \mathcal{L}_{m}
\end{pmatrix}=\begin{pmatrix}\left(z_{11},z_{12},\dots,z_{1N}\right) \\ \left(z_{21},z_{22},\dots,z_{2N}\right) \\ \vdots \\ \left(z_{m1},z_{m2},\dots,z_{mN}\right)  
\end{pmatrix},\hspace{0.1 cm}\mbox{defining the  matrix}\hspace{0.1 cm}Z\overline{Z}^{t}=\begin{pmatrix}\left<\mathcal{L}_{1},\mathcal{L}_{1}\right> & \left<\mathcal{L}_{1},\mathcal{L}_{2}\right>  & \dots  & \left<\mathcal{L}_{1},\mathcal{L}_{m}\right>  \\ \left<\mathcal{L}_{2},\mathcal{L}_{1}\right> &\left<\mathcal{L}_{2},\mathcal{L}_{2}\right> & \dots  & \left<\mathcal{L}_{2},\mathcal{L}_{m}\right>  \\ \vdots & \vdots & \ddots & \vdots \\ \left<\mathcal{L}_{m},\mathcal{L}_{1}\right> & \left<\mathcal{L}_{m},\mathcal{L}_{2}\right>  & \dots  & \left<\mathcal{L}_{m},\mathcal{L}_{m}\right>  \end{pmatrix},\label{row1} \end{equation}

Any Formal (Holomorphic) Embedding, from  $\mathcal{M}_{m,N}$ into $\mathcal{BSD}_{m,N}$, is denoted by
  $\left(G\left(Z,W\right),F\left(Z,W\right)\right)$. Therefore 

\begin{equation}   
  G\left(Z,W\right)= \left<F\left(Z,W\right) ,F\left(Z,W\right) \right>.
   \label{Eq4} \end{equation}

 The equations   (\ref{Eq4})  are used    in order to   implement    linear changes of coordinates preserving the   $\mathcal{BSD}$-Model  (\ref{coordA}). In particular, we use the  product of matrices
\begin{equation} V\otimes Z=\left(\displaystyle\sum_{l=1}^{N}\displaystyle\sum_{k=1}^{m}v_{kl}^{ij}z_{kl}\right)_{1\leq i\leq m\atop 1\leq j\leq N},\quad\mbox{provided a matrix $V=\left( v_{\alpha}^{\beta} \right)_{1\leq\alpha\leq mN}^{1\leq \beta\leq mN}$ defined by (\ref{Ident})  using     the 
identification:}\label{edef}
\end{equation}
\begin{equation*}V\equiv\begin{pmatrix}  {\begin{pmatrix} v^{11}_{11} & v^{11}_{12} & \dots & v^{11}_{1N}\\    v^{12}_{11} & v^{12}_{12} & \dots & v^{12}_{1N} \\  \vdots & \vdots & \ddots & \vdots \\  v^{1N}_{11} & v^{1N}_{12} & \dots & v^{1N}_{1N} \end{pmatrix}} &    \begin{pmatrix}
\dots \\ \dots \\ \ddots \\ \dots
\end{pmatrix}  & {\begin{pmatrix} v^{11}_{m1} & v^{11}_{m2} & \dots & v^{11}_{mN}\\    v^{12}_{m1} & v^{12}_{m2} & \dots & v^{12}_{mN} \\  \vdots & \vdots & \ddots & \vdots \\  v^{1N}_{m1} & v^{1N}_{m2} & \dots & v^{1N}_{mN} \end{pmatrix}} \\   {\begin{pmatrix} \vdots & \vdots & \ddots & \vdots
\end{pmatrix}}& \begin{pmatrix}\vdots\end{pmatrix} & {\begin{pmatrix} \vdots & \vdots & \ddots & \vdots
\end{pmatrix}}    \\ {\begin{pmatrix} v^{m1}_{11} & v^{m1}_{12} & \dots & v^{m1}_{1N}\\    v^{m2}_{11} & v_{12}^{m2} & \dots & v^{m2}_{1N} \\   \dots  & \vdots & \ddots & \vdots \\  v^{mN}_{11} & v^{mN}_{12} & \dots & v^{mN}_{1N} \end{pmatrix}}  & \begin{pmatrix}
\dots \\ \dots \\ \ddots \\ \dots
\end{pmatrix}  & {\begin{pmatrix} v^{m1}_{m1} & v^{m1}_{m2} & \dots & v^{m1}_{mN}\\    v^{m2}_{m1} & v^{m2}_{m2} & \dots & v^{m2}_{mN} \\  \vdots & \vdots & \ddots & \vdots \\  v^{mN}_{m1} & v^{mN}_{m2} & \dots & v^{mN}_{mN} \end{pmatrix}}   \end{pmatrix}\in\mathcal{M}_{mN\times mN}\left(\mathbb{C}\right).
\end{equation*}   

In order to implemented the strategy from \cite{V4},   we reconsider more generally (\ref{edef}).  We consider the linear parts  of   $G\left(0,W\right)$ and 
  of the $F$-component of the Formal Embedding,  denoted by
\begin{equation*}   A=\left({a}^{ij}_{kl}\right)_{1\leq k,l\leq m}^{1\leq i,j\leq m}   \in\mathcal{M}_{m^{2}\times m^{2}}\left(\mathbb{C}\right)\hspace{0.1 cm}\mbox{and}\hspace{0.1 cm}V\equiv \begin{pmatrix} V_{11}  & V_{12}  & \dots &    V_{1\hspace{0.05 cm}m}  \\ V_{21}  & V_{22}  & \dots &    V_{2\hspace{0.05 cm}m}  \\ \vdots & \vdots & \ddots & \vdots  \\ V_{m\hspace{0.05 cm}1}  & V_{m\hspace{0.05 cm}2}  & \dots &   V_{m\hspace{0.05 cm}m} \end{pmatrix}\in\mathcal{M}_{mN\times mN}\left(\mathbb{C}\right).  
\end{equation*} 
 
We make suitable replacements in  (\ref{Eq4}) in order to collect    terms  in $\left(Z,\overline{Z}\right)$. We obtain
 \begin{equation} \quad\hspace{0.25 cm} {a}_{k l }^{i j }\left<\mathcal{L}_{k }, \mathcal{L}_{l }\right>=\left( V_{i k }\left(\mathcal{L}_{k }\right)^{t} \right)\cdot\overline{ \left(  V_{j l }\left(\mathcal{L}_{l }\right)^{t}\right)} , \quad  \mbox{for all $i ,j ,k ,l =1,\dots,m$.}  \label{V33}
  \end{equation}

We  use (\ref{edef}) in the light of (\ref{V33}) according to the proof of Proposition $4.1$ from \cite{V4}.   It follows that up   to compositions with suitable linear   automorphisms of  $\mathcal{BSD}$-Models    (\ref{row1}),    we obtain  
 
 \begin{equation*}    G\left(Z,W\right)= W+\mbox{O}(2) \hspace{0.1 cm}\mbox{and} \hspace{0.1 cm}  F\left(Z,W\right)= \left(
 Z,0\right) +\mbox{O}(2).  
 \end{equation*}

\section{Formal Power Series }
 
  The equation   (\ref{Eq4}) are further studied in order to compute the Formal Embedding    according to the standard linearization procedure (see \cite{V1},\cite{V2},\cite{V4},\cite{CM},\cite{D1},\cite{D3}) in the light of $(2.2)$, $(2.3)$ and $(2.4)$ from the root\cite{V3}. A matrix polynomial (of a certain bidegree) is a matrix whose entries are defined by polynomials of identical bidegree.

  Let $M\subset\mathbb{C}^{mN+m^{2}}$ be a Real-Formal
Submanifold defined near  $p=0$ by
\begin{equation}W=Z\overline{Z}^{t} + \displaystyle\sum_{k+l\geq 3} \varphi_{k,l}\left(Z,\overline{Z}\right),\hspace{0.1 cm}\mbox{provided the homogeneous matrix  polynomials} \label{ew}
\end{equation}
 
\begin{equation*} \varphi_{k,l}\left(Z,\overline{Z}\right):=\begin{pmatrix}  \varphi_{k,l}^{1,1}\left(Z,\overline{Z}\right) &  \dots & \varphi_{k,l}^{1,m}\left(Z,\overline{Z}\right)  \\   \vdots & \ddots & \vdots\\  \varphi_{k,l}^{m,1}\left(Z,\overline{Z}\right) & \dots & \varphi_{k,l}^{m,m}\left(Z,\overline{Z}\right) \end{pmatrix}\hspace{0.1 cm}\mbox{of bidegree  $(k,l)$ in $\left(Z,\overline{Z}\right)$, for all $k,l\in\mathbb{N}$ with $k+l \geq 3$.} \end{equation*}

Provided $(2.6)$ from \cite{V3},   we  consider   the   formal expansions
\begin{equation}F\left(Z,W\right)=\displaystyle\sum_{k,l\geq 0}F_{k,l}\left(Z,W\right)= \begin{pmatrix} \displaystyle\sum_{k,l\geq 0}F_{k,l}^{1,1}\left(Z,W\right)& \dots&\displaystyle\sum_{k,l\geq 0}F_{k,l}^{1,N+1}\left(Z,W\right)   \\  \vdots &   \ddots & \vdots  \\ \displaystyle\sum_{k,l\geq 0}F_{k,l}^{m,1}\left(Z,W\right)& \dots&\displaystyle\sum_{k,l\geq 0}F_{k,l}^{m,N+1}\left(Z,W\right)\end{pmatrix},\label{tz1}\end{equation}

\begin{equation} G\left(Z,W\right)=\displaystyle\sum_{k,l\geq 0}G_{k,l}\left(Z,W\right)=\begin{pmatrix} \displaystyle\sum_{k,l\geq 0}G_{k,l}^{1,1}\left(Z,W\right) &\dots&\displaystyle\sum_{k,l\geq 0}G_{k,l}^{1,m}\left(Z,W\right)   \\    \vdots & \ddots & \vdots  \\ \displaystyle\sum_{k,l\geq 0}G_{k,l}^{m,1}\left(Z,W\right)& \dots&\displaystyle\sum_{k,l\geq 0}G_{k,l}^{m,m}\left(Z,W\right)\end{pmatrix},\label{tz2}\end{equation} such that their entries are homogeneous polynomials of bidegree
$(k,l)$ in $\left(Z,W\right)$,   for all $k,l\in\mathbb{N}$.
 
Combining   (\ref{ew}),(\ref{tz1})  and (\ref{tz2}), we obtain  
 \begin{equation}\begin{split}& \displaystyle\sum_{k,l\geq 0}G_{k,l}\left(Z,Z\overline{Z}^{t} + \displaystyle\sum_{k+l\geq 3} \varphi_{k,l}\left(Z,\overline{Z}\right)\right)=\left(\displaystyle\sum_{k,l\geq 0}F_{k,l}\left(Z,Z\overline{Z}^{t} + \displaystyle\sum_{k+l\geq 3} \varphi_{k,l}\left(Z,\overline{Z}\right)\right)\right)\cdot \hspace{0.2 cm}\\&\quad\quad\quad\quad\quad\quad\quad\quad\quad\quad\quad\quad\quad\quad\quad\quad\quad\quad\quad\overline{\left(\displaystyle\sum_{k,l\geq 0}F_{k,l}\left(Z,Z\overline{Z}^{t} + \displaystyle\sum_{k+l\geq 3} \varphi_{k,l}\left(Z,\overline{Z}\right)\right)\right)}^{t}.\end{split}\label{tx1}
  \end{equation}

Since  the   formal power series  (\ref{tz1}) and (\ref{tz2}) do not have constant terms, we obtain
\begin{equation*}  G_{0,0}\left(Z,W\right)=\mbox{O}_{m\times m} \hspace{0.1 cm}\mbox{and}  \hspace{0.1 cm}   F_{0,0}\left(Z,W\right)=\mbox{O}_{m\times N}. 
 \end{equation*}
  
In order to make further computations in (\ref{tx1}),  we consider the following:

\section{Generalized Fischer-Decompositions}  Motivated by $(2.10)$ from the root\cite{V3}, we recall   the following notation  from Shapiro\cite{Sh}:
\begin{equation}P^{\star}=\displaystyle\sum_{\substack {\left|I_{1}\right|+\left|I_{2}\right|=k_{0}\\ I_{1},I_{2}\in\mathbb{N}^{mN}}}\overline{p}_{I_{1}I_{2}}\frac{\partial^{k_{0}}}{\partial Z^{I_{1}}\partial\overline{Z}^{I_{2}}} ,\quad\mbox{provided $P\left(Z,\overline{Z}\right)=\displaystyle\sum_{\substack {\left|I_{1}\right|+\left|I_{2}\right|=k_{0}\\I_{1},I_{2}\in\mathbb{N}^{mN}}}p_{I_{1}I_{2}}Z^{I_{1}}\overline{Z}^{I_{2}}$, where $k_{0}\in\mathbb{N}$.} \label{Pstar}\end{equation}

Provided   the space of all homogeneous polynomials of degree $k$ in $Z$ denoted by $\mathbb{H}_{k}$, we recall   the Fischer inner product from Shapiro\cite{Sh}. In particular,   $(2.11)$ from the root\cite{V3} must be reformulated like
\begin{equation}\left<Z^{I_{1}};\hspace{0.1 cm}Z^{I_{2}}\right>_{\mathcal{F}}=\left\{\substack{\hspace{0.05 cm} 0,\quad\hspace{0.15 cm} I_{1}\neq I_{2} \\ I_{1}!,\quad I_{1}=I_{2}}\right.,\quad
\mbox{for all  $I_{1},I_{2}\in\mathbb{N}^{mN}$. }\label{ppo}\end{equation}

 In particular,  we  reformulate the first part of Lemma $2.1$ from \cite{V4}: we define $$\mathcal{J}_{p}=\left\{J\in\mathbb{N}^{mN};\hspace{0.1 cm} \left|J\right|=p\right\}.$$

\bl \label{L2.1} Let
$P\left(Z,\overline{Z}\right)$ be a bihomogeneous polynomial of bidegree $(p,q)$ in  $\left(Z,\overline{Z}\right)$ with $p> q$.  Uniquely there exist the polynomials 
$$\mbox{$\left\{Q_{J}\left(  Z \right)\right\}_{J\in\mathcal{J}_{q}}$
 and $R\left(Z,\overline{Z}\right)\in   \displaystyle \bigcap_{J\in \mathcal{J}_{p}}\ker\left(  \displaystyle\prod_{k,l=1}^{m}\langle
\mathcal{L}_{k},\mathcal{L}_{l}\rangle^{j_{kl}}\right)^{\star} $, such that:} $$ 
  
\begin{equation*}  P\left(Z,\overline{Z}\right)=\displaystyle\sum_{J\in \mathcal{J}_{q}}Q_{J}\left(Z\right)\cdot\displaystyle\prod_{k,l=1}^{m}\langle
\mathcal{L}_{k},\mathcal{L}_{l}\rangle^{j_{kl}} +R\left(Z,\overline{Z}\right). \end{equation*}  \el 
 
Before going forward, we  reformulate the second part Lemma $2.1$ from \cite{V4}:
\bl \label{L2.2}
Let
$P\left(Z,\overline{Z}\right)$ be a bihomogeneous polynomial of bidegree $(p,q)$ in  $\left(Z,\overline{Z}\right)$ with $p \leq q$. Uniquely there exist the  polynomials
$$\mbox{$\left\{Q_{J}^{(ij)}\left(Z\right)\right\}_{J\in\mathcal{J}_{p-1}\atop{i=1,\dots,m\atop{j=1,\dots,N}}}$ and $R\left(Z,\overline{Z}\right)\in \displaystyle \bigcap_{i=1}^{m} \bigcap_{j=1}^{N}  \displaystyle \bigcap_{J\in \mathcal{J}_{p-1}}\ker\left(z_{ij} \cdot \displaystyle\prod_{k,l=1}^{m}\langle
\mathcal{L}_{k},\mathcal{L}_{l}\rangle^{j_{kl}}\right)^{\star}$,} 
\hspace{0.1 cm}\mbox{such that:}$$
 
\begin{equation*} P\left(Z,\overline{Z}\right)=\displaystyle\sum_{i=1}^{m}\displaystyle\sum_{j=1}^{N}z_{ij}\displaystyle\sum_{J\in \mathcal{J}_{p-1}}\overline{\left(Q^{(ij)}_{I}\left(Z\right)\right)}\cdot\displaystyle\prod_{k,l=1}^{m}\langle
\mathcal{L}_{k},\mathcal{L}_{l}\rangle^{j_{kl}}  +R\left(Z,\overline{Z}\right). \end{equation*}

\el
 
 Their proofs follows from the proof of Lemma 2.1 in \cite{V3}. We are ready to recall the  Fischer norm \cite{Sh}  defined by
\begin{equation}\left\|f_{k}\left(z\right)\right\|:=\displaystyle\sum_{\left|I\right|=k\atop{I\in\mathbb{N}^{N}}}I!\left|c_{I}\right|^{2},\quad\mbox{if  $f_{k}\left(z\right):=\displaystyle\sum_{\left|I\right|=k\atop{I\in\mathbb{N}^{N}}}c_{I}z^{I}$.}\label{fnorm}\end{equation}

Provided $f\left(Z\right)$, $g\left(Z\right)\in\mathbb{H}_{k}$ defining the orthogonal decomposition $f\left(Z\right)=g\left(Z\right)+h\left(Z\right)$, we recall  from Shapiro\cite{Sh} that $$\left\|f\left(Z\right)\right\|=\left\|g\left(Z\right)\right\|+\left\|h\left(Z\right)\right\|.$$

  \section{Computations}
  
We pursue computations of normal form type motivated by   the partial normal form  \cite{V1}.
 The chosen  Fischer Decompositions are motivated by how the formal transformation   appears in the equation (\ref{tx1}). Other strategies may be possibly considered as in  \cite{V3}.

In order to make computations, we write formal power series expansion 
\begin{equation}\left(Z',W'\right)=\left(Z+\displaystyle\sum_{p+q\geq 2}F_{p,q}\left(Z,W\right),W+\displaystyle\sum_{p+q
\geq2}G_{p,q}\left(Z,W\right)\right),\label{EMLmap}\end{equation} where
$F_{p,q}\left(Z,W\right)$, $G_{p,q}\left(Z,W\right)$ are homogeneous polynomials in $Z$ of
degree $p$ and degree $q$ in $W$.

We compute   the polynomials $F_{m,n'}\left(Z\right)$ with $m+2n'=T-1$, and respectively
$G_{k,l}\left(Z\right)$ with $k+2l=T$ using  induction   on $T \geq 3$. In particular, we assume  that we have computed the polynomials $F_{k,l}\left(Z\right)$ with
$k+2l<T-1$, $G_{k,l}\left(Z\right)$ with $k+2l<T$. We collect the terms of bidegree $(p,q)$ in $ \left(Z,\overline{Z}\right)$
with $T=p+q$ in (\ref{tx1}) in analogy to $(2.21)$ from \cite{V3}. We  obtain   

\begin{equation}\begin{split}&
-\begin{pmatrix}  \varphi _{p,q}^{1,1} &   \dots & \varphi ^{1,m}_{p,q}    \\   \vdots & \ddots & \vdots \\ \varphi _{p,q}^{m,1} &  \dots & \varphi ^{m,m}_{p,q} \end{pmatrix}\left(Z,\overline{Z}\right)=\displaystyle\sum_{  I \in\mathcal{I}_{q}}\begin{pmatrix} G^{1,1}_{p-q,I}\left(Z\right) &  \dots & G^{1,m}_{p-q,I}\left(Z\right)  \\    \vdots &   \ddots & \vdots\\  G^{m,1}_{p-q,I}\left(Z\right) & \dots & G^{m,m}_{p-q,I}\left(Z\right)\end{pmatrix}\cdot\displaystyle\prod_{k,l=1}^{m}\langle
\mathcal{L}_{k},\mathcal{L}_{l}\rangle^{i_{kl}}\\&\quad\quad\quad\quad\quad -\displaystyle\sum_{I \in\mathcal{I}_{q-1}}\begin{pmatrix} \left\langle F^{1}_{p-q+1,I}\left(Z\right),\mathcal{L}_{1}\right\rangle   & \dots &\left\langle F^{1}_{p-q+1,I}\left(Z\right),\mathcal{L}_{m}\right\rangle     \\ \vdots &  \ddots & \vdots \\ \left\langle F^{m}_{p-q+1,I}\left(Z\right),\mathcal{L}_{1}\right\rangle   & \dots &\left\langle F^{m}_{p-q+1,I}\left(Z\right),\mathcal{L}_{m}\right\rangle \end{pmatrix}\cdot\displaystyle\prod_{k,l=1}^{m}\langle
\mathcal{L}_{k},\mathcal{L}_{l}\rangle^{j_{kl}}\\&\quad\quad\quad\quad\quad -\displaystyle\sum_{I \in\mathcal{I}_{q-1}}\begin{pmatrix} \overline{\left\langle F^{1}_{m-n+1,I}\left(Z\right),\mathcal{L}_{1}\right\rangle} &  \dots & \overline{\left\langle F^{m}_{p-q+1,I}\left(Z\right),\mathcal{L}_{1}\right\rangle} \\   \vdots &   \ddots & \vdots  \\  \overline{\left\langle F^{1}_{p-q+1,I}\left(Z\right), \mathcal{L}_{m}\right\rangle} &   \dots & \overline{\left\langle F^{m}_{p-q+1,I}\left(Z\right),\mathcal{L}_{m}\right\rangle}\end{pmatrix}\cdot\overline{\displaystyle\prod_{k,l=1}^{m}\langle
\mathcal{L}_{k},\mathcal{L}_{l}\rangle^{i_{kl}} }.\end{split}
\label{ecg3}
\end{equation}
 
In analogy to $(2.22)$ from \cite{V3} while $p<q-1$,  we  obtain 
\begin{equation}  0=-\displaystyle\sum_{J \in\mathcal{I}_{q-1}}\left\langle
\mathcal{L}_{1} ,F^{1}_{p-q+1,J}(Z)\right\rangle\cdot\displaystyle\prod_{k,l=1}^{m}\langle
\mathcal{L}_{k},\mathcal{L}_{l}\rangle^{j_{kl}} +\displaystyle\sum_{i=1}^{m} \left(\left(\varphi_{p,q}\right)^{i,i}\right)\left(Z,\overline{Z}\right). \label{11121} \end{equation}

According to Lemma \ref{L2.2}, write the Fischer Decomposition
\begin{equation}\displaystyle\sum_{i=1}^{m} \left( \left(\varphi_{p,q}\right)^{i,i}\right)\left(Z,\overline{Z}\right)= -\displaystyle\sum_{i=1}^{m}\displaystyle\sum_{j=1}^{N}z_{ij} \left(
\displaystyle\sum_{J \in\mathcal{I}_{q-1}}\overline{Q^{j}_{J}\left(Z\right)} \cdot\displaystyle\prod_{k,l=1}^{m}\langle
\mathcal{L}_{k},\mathcal{L}_{l}\rangle^{j_{kl}}\right) +R'\left(Z,\overline{Z}\right),\quad\mbox{such that:}\label{1111}  \end{equation}

$$
R'\left(Z,\overline{Z}\right)\in \displaystyle \bigcap_{i=1}^{m}  \displaystyle \bigcap_{j=1}^{N}\displaystyle \bigcap_{J \in\mathcal{I}_{q-1}}\ker\left( z_{1j} \displaystyle\prod_{k,l=1}^{m}\langle
\mathcal{L}_{k},\mathcal{L}_{l}\rangle^{j_{kl}}\right)^{\star}.$$

The uniqueness of the Fischer Decomposition yields
\begin{equation}F^{1}_{p-q+1,J}\left(Z,W\right)=\displaystyle\sum_{J \in\mathcal{I}_{q-1}}\left( Q^{1}_{J},\dots,Q^{N}_{J}\right)\left(Z\right)W^{J} .\label{F1}\end{equation}

In analogy to $(2.25)$ from \cite{V3},  we  obtain 
 
\begin{equation}  0=-\displaystyle\sum_{J \in\mathcal{I}_{q-1}}\left\langle
\mathcal{L}_{1} ,F^{1}_{2,J}(Z)\right\rangle\cdot\displaystyle\prod_{k,l=1}^{m}\langle
\mathcal{L}_{k},\mathcal{L}_{l}\rangle^{j_{kl}} -\displaystyle\sum_{J \in\mathcal{I}_{q-1}}\left\langle
F^{1}_{0,J}(Z),z^{1}\right\rangle\cdot\displaystyle\prod_{k,l=1}^{m}\langle
\mathcal{L}_{k},\mathcal{L}_{l}\rangle^{j_{kl}}+\displaystyle\sum_{i=1}^{m}\left( \left(\varphi_{q,q+1}\right)^{i,i}\right)\left(Z,\overline{Z}\right).\label{1112}  \end{equation}

According to Lemma \ref{L2.2}, write the Fischer Decomposition

\begin{equation} \displaystyle\sum_{i=1}^{m} \left( \left(\varphi_{q,q+1}\right)^{i,i}\right)\left(Z,\overline{Z}\right)= \displaystyle\sum_{i=1}^{m}\displaystyle\sum_{j=1}^{N}
 z_{ij}  \left(\displaystyle\sum_{J \in\mathcal{I}_{q-1}}
\overline{C^{j}_{J}\left(Z\right)}\cdot\displaystyle\prod_{k,l=1}^{m}\langle
\mathcal{L}_{k},\mathcal{L}_{l}\rangle^{j_{kl}}\right) +R_{1}'\left(Z,\overline{Z}\right),\quad\mbox{such that:}\label{111}  \end{equation}
 
 $$
R_{1}'\left(Z,\overline{Z}\right)\in   \displaystyle \bigcap_{i=1}^{m} \displaystyle \bigcap_{j=1}^{N}\displaystyle \bigcap_{J \in\mathcal{I}_{q-1}}\ker\left( z_{ij} \displaystyle\prod_{k,l=1}^{m}\langle
\mathcal{L}_{k},\mathcal{L}_{l}\rangle^{j_{kl}}\right)^{\star}.$$

The uniqueness of the Fischer Decomposition yields
\begin{equation}\left(F^{1}_{2,q-1}\left(Z,W\right),\dots,F^{m}_{2,q-1}\left(Z,W\right)\right)=\displaystyle\sum_{J \in\mathcal{I}_{q-1}}\left(C^{11}_{J},\dots,C^{1N}_{J},\dots,C^{m1}_{J},\dots,C^{mN}_{J}\right)\left(Z\right)W^{J}.\label{F2}\end{equation}

   Collecting the terms of bidegree $(p,q)$ in
$\left(Z,\overline{Z}\right)$ in (\ref{ecg3}) with $p\geq q$ and $p,q\geq 1$, we
obtain   
\begin{equation} -
\varphi_{p,q}\left(Z,\overline{Z}\right)=\displaystyle\sum_{J \in\mathcal{I}_{q}} G_{p-q,J}\left(Z\right)\cdot\displaystyle\prod_{k,l=1}^{m}\langle
\mathcal{L}_{k},\mathcal{L}_{l}\rangle^{j_{kl}}-  \displaystyle\sum_{I \in\mathcal{I}_{q-1}}\begin{pmatrix} \left\langle F^{1}_{p-q+1,I}\left(Z\right),\mathcal{L}_{1}\right\rangle & \dots & \left\langle F^{1}_{p-q+1,I}\left(Z\right),\mathcal{L}_{m}\right\rangle    \\ \vdots &  \ddots & \vdots  \\ \left\langle F^{m}_{p-q+1,I}\left(Z\right),\mathcal{L}_{1}\right\rangle &   \dots & \left\langle F^{m}_{p-q+1,I}\left(Z\right),\mathcal{L}_{m}\right\rangle  \end{pmatrix}\cdot\displaystyle\prod_{k,l=1}^{m}\langle
\mathcal{L}_{k},\mathcal{L}_{l}\rangle^{i_{kl}}.\label{1A1} 
\end{equation}

According to Lemma \ref{L2.1}, write the Fischer Decomposition
 \begin{equation}\begin{split}&\varphi_{p,q}\left(Z,\overline{Z}\right)+ \displaystyle\sum_{J \in\mathcal{I}_{q-1}}\begin{pmatrix} \left\langle F^{1}_{p-q+1,I}\left(Z\right),\mathcal{L}_{1}\right\rangle &  \dots & \left\langle F^{1}_{p-q+1,I}\left(Z\right),\mathcal{L}_{m}\right\rangle   \\ \vdots& \ddots & \vdots  \\ \left\langle F^{m}_{p-q+1,I}\left(Z\right),\mathcal{L}_{1}\right\rangle  & \dots & \left\langle F^{m}_{p-q+1,I}\left(Z\right),\mathcal{L}_{m}\right\rangle  \end{pmatrix} \cdot\displaystyle\prod_{k,l=1}^{m}\langle
\mathcal{L}_{k},\mathcal{L}_{l}\rangle^{j_{kl}}\\&\quad = \displaystyle\sum_{J \in\mathcal{I}_{q}} E_{m,J}\left(Z\right)\cdot\displaystyle\prod_{k,l=1}^{m}\langle
\mathcal{L}_{k},\mathcal{L}_{l}\rangle^{j_{kl}}+\begin{pmatrix} R_{11}'\left(Z,\overline{Z}\right) &   \dots & R_{1m}'\left(Z,\overline{Z}\right)  \\ \vdots   & \ddots & \vdots  \\ R_{m1}'\left(Z,\overline{Z}\right)&   \dots & R_{mm}'\left(Z,\overline{Z}\right) \end{pmatrix},\quad\mbox{such that:}\end{split}\label{1A}\end{equation}

$$
R_{11}'\left(Z,\overline{Z}\right),R_{12}'\left(Z,\overline{Z}\right),\dots\in  \displaystyle \bigcap_{J \in\mathcal{I}_{q}}\ker\left(\displaystyle\prod_{k,l=1}^{m}\langle
\mathcal{L}_{k},\mathcal{L}_{l}\rangle^{j_{kl}}\right)^{\star}.$$

The uniqueness of the Fischer Decomposition yields
\begin{equation}G_{p-q,q}\left(Z,W\right)=\displaystyle\sum_{J \in\mathcal{I}_{q}} E_{m,J}\left(Z\right)W^{J}.\label{gig}\end{equation}

  Collecting  the terms of bidegree $\left(T,0\right)$ and
$\left(0,T\right)$ in $\left(Z,\overline{Z}\right)$ in  (\ref{ecg3}), we
obtain  
     
\begin{equation}G_{T,0}\left(Z\right)=\varphi_{T,0}\left(Z\right)-\left(\varphi_{0,T}\left(  Z\right)\right)^{t}.\label{F3}\end{equation}

\section{Further Normalizations}

We consider compositions with suitable automorphisms of the  $\mathcal{BSD}$-Model $\mathcal{M}$ from (\ref{coordA}) in order to assume
\begin{equation}
F_{0,q+1}\left(Z,W\right)=0 \hspace{0.1 cm}\mbox{and}  \hspace{0.1 cm} F_{1,q}\left(Z,W\right)=0,\quad \mbox{for all}\hspace{0.1
cm} q\geq 1. \label{fn}
\end{equation} 
  
 It is not possible to adapt the  formula of  automorphisms  computed by Huang-Yin\cite{HuYi} in $\mathbb{C}^{N+1}$, but it is possible to continue these computations in order to extended the result even in the non-equidimensional case described by Theorem \ref{teo1}. Other details to be added.

\section{Rapid Iterations}

 \subsection{Notations}   We consider $R:=\left(r ,\dots,r \right)$ in order to   define 
\begin{equation*} \Delta_{r}=\left\{\left(Z,W\right)\in\mathbb{C}^{mN+m^{2}};\hspace{0.1 cm}\left|z_{ij}\right|<r \hspace{0.1 cm}\mbox{and}\hspace{0.1 cm} \left|w_{kl}\right|^{2}<Nr^{2},\hspace{0.1 cm}\mbox{for all $k,l,i=1,\dots,m$ and $j=1,\dots,N$.}\right\} \end{equation*}

Provided $\xi$   defined similarly as $Z$,  we define
\begin{equation*}        D_{r}=\left\{\left(Z,\xi\right)\in\mathbb{C}^{2N}\times\mathbb{C}^{2N},\hspace{0.1 cm}\left|z_{ij}\right|<r ,\hspace{0.1 cm}\left|\xi_{ij}\right|<r,\hspace{0.1 cm}\mbox{for all $\left(i,j\right)\in\left\{1,\dots,N\right\}\times\left\{1,2\right\}$}\right\}.\end{equation*}
 
Provided $E\left(Z,\xi\right)$ a holomorphic function defined over $ \overline{D}_{r}$, and  $h\left(Z,W\right)$  a holomorphic function defined over $\overline{\Delta}_{r}$, we define
 \begin{equation*}\left\|E\right\|_{r}:=\displaystyle\sup_{\left(Z,\xi\right)\in  \overline{D}_{r}}\left|E\left(Z,\xi\right)\right|\hspace{0.1 cm}\mbox{and}\hspace{0.1 cm}\left|h\right|_{r}=\displaystyle\sup_{\left(Z,W\right)\in \Delta_{r}}\left|h\left(Z,W\right)\right|.\end{equation*}

More generally, provided the matrix $E\left(Z,\overline{Z}\right)$  defined by  \begin{equation*}E\left(Z,\overline{Z}\right)=\left(E_{ij}\left(Z,\overline{Z}\right) \right)_{1\leq i,j \leq m},\quad\mbox{we use the notation}\hspace{0.1 cm}\left\|E\right\|_{r}=\max_{1\leq i,j \leq m}\left\|E_{ij}\right\|_{r}.\end{equation*}
 
\subsection{ Uploads from Moser\cite{Mo}} We  consider the numbers
\begin{equation}\frac{1}{2}<r'<\sigma<\rho<r\leq 1,\quad \rho =\frac{2 r'+r}{3},\quad \sigma =\frac{2r'+\rho}{3},\quad n\in\mathbb{N}^{\star}.\label{num}\end{equation}

In particular, we have
\begin{equation}\left(\frac{\tau}{r}\right)^{2}\leq\frac{\rho}{r},\quad\mbox{for all $\frac{1}{2}<\rho<\tau<r \leq1$ and $\tau=\frac{r+2\rho}{3}$}.\label{xw}\end{equation}

In particular, we define the following sequences of numbers
  $$r_{n}:=\frac{1}{2}\left(1+\frac{1}{n+1}\right),\quad \rho_{n}=\frac{r_{n+1}+2r_{n}}{3},\quad 
 \sigma_{n}=\frac{\rho_{n}+2r_{n}}{3}.$$

In particular, we assume $r=r_{n}$, $\rho=\rho_{n}$, $r'=r_{n+1}$, for all $n\in\mathbb{N}$.  We observe
\begin{equation} \frac{r_{n+1}}{r_{n}}=1-\frac{1}{\left(n+1\right)^{2}},\quad \frac{1}{r_{n}-r_{n+1}}=\left(n+1\right)\left(n+2\right).\label{rrr}\end{equation}

 \subsection{Estimations} We consider the Real-Analytic Submanifold $M\subset\mathbb{C}^{mN+m^{2}}$ defined near $p=0$ by
\begin{equation}W=\Phi\left(Z,\overline{Z}\right)=Z\overline{Z}^{t}+E\left(Z,\overline{Z}\right),\label{r1}\end{equation}

  Following the strategy of Huang-Yin\cite{HuYi},  we define  
\begin{equation}    \Theta\left(Z,W\right):=\left(Z+\widehat{F}\left(Z,W\right),\hspace{0.1 cm}W+\widehat{G}\left(Z,W\right)\right)=\left(Z+F_{\rm{nor}}^{d-1}\left(Z,W\right)+O_{\rm{wt}}\left(d\right),\hspace{0.1 cm}W+G_{\rm{nor}}^{\left(d\right)}\left(Z,W\right)+\rm{O}_{\rm{wt}}\left(d+1\right)\right), \label{r2}\end{equation}
which sends $M$ up to the degree $d$ into the model manifold $M_{\infty}$ defined by (\ref{15}):
\begin{equation}M'=\Theta\left(M\right):\hspace{0.1 cm} W'=Z'{\overline{Z'}}^{t}.\label{r3}\end{equation}

Assume  $\rm{Ord}\left(E\left(Z,\xi\right)\right)\geq d$. If $\Theta$ is defined in (\ref{r2}) and $m$ is defined in (\ref{r3}), then $\rm{Ord}\left(E'\left(Z,\xi\right)\right)\geq 2d-2$, where
\begin{equation}M': W'=Z'{\overline{Z'}}^{t}+E'\left(Z,\overline{Z}\right).\label{900}\end{equation}

Let $J^{2d-3}\left(E\left(Z,\overline{Z}\right)\right)$ be  the polynomial defined by the Taylor expansion of $E\left(Z,\overline{Z}\right)$ up to the degree $2d-3$. Therefore
\begin{equation} \begin{split}&E'\left(Z',\overline{Z'}\right)=G\left(Z,\Phi\left(Z,\overline{Z}\right)\right)-G\left(Z,W_{0}\right)-2\Re\left\{Z \overline{\left(F\left(Z,\Phi\left(Z,\overline{Z}\right)\right)-F\left(Z,W_{0}\right)\right)}^{t} \right\}+
\\&\quad\quad\quad\quad\quad\quad\quad\left(F\left(Z,\Phi\left(Z,\overline{Z}\right)\right)\right)\left(\overline{F\left(Z,\Phi\left(Z,\overline{Z}\right)\right)}\right)^{t}-\left(J^{2d-3}\left(E\left(Z,\overline{Z}\right)\right)-
E\left(Z,\overline{Z}\right)\right),
 \end{split}\label{r4}\end{equation}

\begin{equation}
\left\|E -J^{2d-3}\left(E \right)\right\|_{\rho}
\leq\frac{\left(2d\right)^{4N}\left\|E\right\|_{r}}{\left(r-\rho\right)^{2N}}\left(\frac{\rho}{r}\right)^{2d-2}.
\end{equation}

\begin{equation} \left|F_{k,l}\left(Z,W\right)\right|_{\rho}\leq  \frac{4}{N}\left(2d\right)^{4N}\left\|E\right\|_{r} \left(\frac{\rho}{r}\right)^{2d-3},\quad\quad
 \left|\widehat{G}_{\alpha,\beta}
\left(Z,W\right)\right|_{\rho}\leq   \left(\left(2d\right)^{4N}+\left(2d\right)^{6N}\right)\left\|E\right\|_{r} \left(\frac{\rho}{r}\right)^{2d-2} ,\label{1003}\end{equation}
$\mbox{for all $k\in\left\{1,\dots,N\right\}$}$, $\alpha,\beta,l\in\{1,2\}$.

 \begin{proof} We obtain
 \begin{equation}\left\|E\left(Z,\xi\right)-J^{2d-3}\left(E\left(Z,\xi\right)\right)\right\|_{\rho}\leq \left\|\displaystyle\sum_{\substack{\left|I\right|+\left|J\right|\geq 2d+2\\ I,J\in\mathbb{N}^{2N}}}a_{I,J}Z^{I}\overline{Z}^{J}\right\|_{\rho}\leq \displaystyle\sum_{\substack{\left|I\right|+\left|J\right|\geq 2d+2\\ I,J\in\mathbb{N}^{2N}}}\left\|E\right\|_{r}\left(\frac{R'}{R}\right)^{I+J} \leq
 \frac{\left(2d\right)^{4N}\left\|E\right\|_{r}}{\left(r-\rho\right)^{2n}}\left(\frac{\rho}{r}\right)^{2d-2},\end{equation}
where  $R':=\left( \rho,\dots, \rho\right)$ and $R:=\left( r,\dots, r\right) $.

By (\ref{11121}) together with (\ref{1111}), (\ref{F1}) and (\ref{1112}), (\ref{F2})  we obtain the following
\begin{equation}\left|\widehat{F}_{k,l}\left(Z,W\right)\right|_{\rho}\leq \frac{4\left(2d\right)^{4N}\left\|E\right\|_{r}}{N}   \left(\frac{\rho}{r}\right)^{2d-3},\quad\mbox{for all $k=1,\dots,N$ and $l\in\{1,2\}$.}\label{s222}\end{equation}

The main  ingredient for computing the $G$-component of our transformation is the following remark
\br Let $S\left(Z,\overline{Z}\right)$ be a homogeneous polynomial of degree $k$ in $\left(Z,\overline{Z}\right)$ written as follows
$$S\left(Z,\overline{Z}\right)=\displaystyle\sum_{\substack{\left|I\right|+\left|J\right|=k\\ I,J\in\mathbb{N}^{2N}}}c_{I,J}Z^{I}\overline{Z}^{J}.$$
Then the following holds
$$\left|S\left(Z,\overline{Z}\right)\right|^{2}\leq\frac{\left\|S\left(Z,\overline{Z}\right)\right\|^{2}}{k!}
\left(\left|z_{11}\right|^{2}+\dots+\left|z_{1N}\right|^{2}+\left|z_{21}\right|^{2}+\dots+\left|z_{2N}\right|^{2}\right)^{2k},$$
and as well the following Cauchy estimates using the domain $D_{r}$ 
$$\left\|S\left(Z,\overline{Z}\right)\right\|^{2}_{r}\leq \frac{k! (k+1)^{2N}}{r^{2k}}\|S\|^{2}_{r}. $$
\er
\begin{proof} From the Cauchy inequality  we obtain the first inequality. By applying the Cauchy formulas using the domain $D_{r}$ following Shapiro\cite{Sh} ,  we obtain 
$$\left\|S\left(Z,\overline{Z}\right)\right\|^{2}=\displaystyle\sum_{\substack{\left|I\right|+\left|J\right|=k\\ I,J\in\mathbb{N}^{2N}}}I!J!\left|c_{I,J}\right|^{2}\leq \frac{\|S\|^{2}_{r}}{r^{2k}}\left(
\displaystyle\sum_{\substack{\left|I\right|+\left|J\right|=k\\ I,J\in\mathbb{N}^{2N}}}I!J!\right)\leq \frac{k! (k+1)^{2N}}{r^{2k}}\left\|S\right\|^{2}_{r}.$$
\end{proof}

 By (\ref{1A1}) using the previous remark together with Lemma $2.2$,  we obtain 
\begin{equation}\left|\widehat{G}_{\alpha,\beta}\left(Z,W\right)\right|_{\rho}\leq   \left(2d\right)^{4N}\left(1  + \left(2d\right)^{2N}\right)\left\|E\right\|_{r} \left(\frac{\rho}{r}\right)^{2d-2},\quad\mbox{for all  $\alpha,\beta\in\{1,2\}$,}\end{equation}
 
\end{proof}

  $d_{n}:=\rm{Ord}\left(E_{n}\left(Z,\overline{Z}\right)\right)\geq 2^{n}+2$, for all $n\in\mathbb{N}$.  

  \section{Acknowledgements}  I was very unbalanced purposefully prior to the  submission of \cite{V3} in its final form.  I wrote this new version motivated by the unbalances from \cite{V3}. My paper \cite{V3} is essentially a failure,  but a very inspirational failure. This work continuation is a very important engine for my research-program, because  I have been able recently to better understand the rapid iteration procedure of Moser. Such procedure may be crucial even in the non-equidimensional case. This failure helped me to understand the rapid iterations much better. I will add more details when I will be calm. To be submitted in its final form at Sinica, because I want them to see what happened to me. The reparation of my academic career remains a very serious open problem.

\end{document}